\documentclass[11pt]{amsart}

\usepackage{amsmath,amssymb,color}

\usepackage[all]{xy}

\newcommand{\C}{{\mathbb{C}}}

\newcommand{\N}{{\mathbb{N}}}

\newcommand{\Ch}{{\mathcal C}}
\newcommand{\Dh}{{\mathcal D}}

\newcommand{\Fh}{{\mathcal F}}
\newcommand{\Gh}{{\mathcal G}}

\newcommand{\Oh}{{\mathcal O}}

\newcommand{\Qh}{{\mathcal Q}}

\newcommand{\Uh}{{\mathcal U}}

\newcommand{\Zh}{{\mathcal Z}}

\newcommand{\au}{\approx_{\mathrm{u}}}
\newcommand{\asu}{\approx_{\mathrm{uh}}}
\newcommand{\sau}{\approx_{\mathrm{su}}}
\newcommand{\sasu}{\approx_{\mathrm{suh}}}
\newcommand{\Aut}{\mathrm{Aut}\,}
\newcommand{\be}{\mathbf{1}}

\newcommand{\halb}{\frac{1}{2}}

\newcommand{\id}{\mathrm{id}}

\newcommand{\ot}{\otimes}

\newcounter{number}[section]

\newenvironment{nummer}{\refstepcounter{number}
{\noindent\arabic{section}.\arabic{number}}}{}

\newcommand{\bn}{\begin{nummer} \rm}
\newcommand{\en}{\end{nummer}}

\newenvironment{thms}{\noindent {\sc Theorem:} \it}{}
\newenvironment{lms}{\noindent {\sc Lemma:} \it}{}
\newenvironment{props}{\noindent {\sc Proposition:} \it}{}
\newenvironment{dfs}{\noindent {\sc Definition:} \it}{}
\newenvironment{cors}{\noindent {\sc Corollary:} \it}{}

\newenvironment{rems}{\noindent {\sc Remark:}}{}

\newenvironment{nproof}{\noindent {\sc Proof:}}{\mbox{}\hfill
\rule[-.2ex]{.25em}{1.8ex}}

\begin{document}

\title{\sc On the $KK$-theory of strongly self-absorbing  $C^{*}$-algebras}

\author{Marius Dadarlat}
\address{Department of Mathematics, Purdue University,
West Lafayette,\\
IN 47907, USA}

\email{mdd@math.purdue.edu}

\author{Wilhelm Winter}
\address{Mathematisches Institut der Universit\"at M\"unster\\
Einsteinstr.\ 62\\ D-48149 M\"unster, \indent Germany}

\email{wwinter@math.uni-muenster.de}

\date{\today}
\thanks{{\it Supported by:} The first named author was partially supported by
NSF grant \#DMS-0500693. \\
\indent The second named author was supported by the DFG (SFB 478).}

\setcounter{section}{-1}

\begin{abstract} 
Let $\Dh$ and $A$ be  unital and  separable $C^{*}$-algebras; let $\Dh$ be 
 strongly self-absorbing. It is known that any two unital $^*$-homomorphisms  from $\Dh$ to $A \otimes \Dh$ are approximately unitarily equivalent. We show that, if $\Dh$ is also $K_{1}$-injective, they are even asymptotically unitarily equivalent. This in particular implies that any unital endomorphism of $\Dh$ is asymptotically inner. Moreover,  the space of automorphisms of $\Dh$ is compactly-contractible (in the point-norm topology) in the sense that
  for any compact Hausdorff space $X$, the set of homotopy classes
  $[X,\Aut(\Dh)]$ reduces to a point. The respective statement 
  holds for the space of unital endomorphisms of $\Dh$.
 As an application,  we give a description of the Kasparov group $KK(\Dh, A\ot \Dh)$
  in terms of $^*$-homomorphisms and  asymptotic unitary equivalence.  
Along the way, we
 show that the Kasparov group $KK(\Dh, A\ot \Dh)$ is isomorphic to $K_0(A\ot \Dh)$.
\end{abstract}

\maketitle

\section{Introduction}

A unital and separable $C^*$-algebra $\Dh \neq \C$ is strongly self-absorbing if there is an isomorphism $\Dh\stackrel\sim\to\Dh\ot \Dh$
 which is approximately unitarily equivalent
to the inclusion map $\Dh\to \Dh\ot \Dh$, $d\mapsto d \ot \be_{\Dh}$ (\cite{TomsWinter:ssa}).
 Strongly self-absorbing $C^{*}$-algebras are known to be simple and nuclear; moreover, they are either purely infinite or stably finite.
The only known examples of strongly self-absorbing $C^{*}$-algebras are the UHF algebras of infinite type (i.e., every prime number that occurs in the respective supernatural number occurs with infinite multiplicity), the Cuntz algebras $\Oh_{2}$ and $\Oh_{\infty}$, the Jiang--Su algebra $\Zh$ and tensor products of $\Oh_{\infty}$ with UHF algebras of infinite type, see \cite{TomsWinter:ssa}. All these examples are  $K_{1}$-injective,
 i.e., the canonical map $\Uh(\Dh)/\Uh_{0}(\Dh) \to K_{1}(\Dh)$ is injective.

It was observed in \cite{TomsWinter:ssa} that any two unital $^*$-homomorphisms $\sigma,\gamma: \Dh \to A \otimes \Dh$ are approximately unitarily equivalent, were $A$ is another unital and separable $C^{*}$-algebra. If $\Dh$ is $K_{1}$-injective, the unitaries implementing the equivalence may even be chosen to be homotopic to the unit. When $\Dh$ is $\Oh_{2}$, $\Oh_{\infty}$, it was known that $\sigma$ and $\gamma$ are even asymptotically unitarily equivalent -- i.e., they can be intertwined by a continuous path of unitaries, parametrized by a half-open interval. Up to this point, it was not clear whether the respective statement holds for the Jiang--Su algebra $\Zh$. Theorem~\ref{D-asu} below provides an affirmative answer to this problem. Even more, we show that the path intertwining $\sigma$ and $\gamma$ may be chosen in the component of the unit. 

We believe this result, albeit technical, is interesting in its own right, and that it will be a useful ingredient for the systematic further use of strongly self-absorbing $C^{*}$-algebras in Elliott's program to classify nuclear $C^{*}$-algebras by $K$-theory data. In fact, this point of view is our main motivation for the study of strongly self-absorbing $C^{*}$-algebras; see \cite{Kir:class}, \cite{Phi:class}, \cite{Winter:Z-stableTAF}, \cite{Winter:lfdr}, \cite{Winter:localizingEC} and \cite{TomsWinter:ASH}    for already existing results in this direction.

For the time being, we use Theorem~\ref{D-asu} to derive some consequences for the Kasparov groups of the form $KK(\Dh, A \otimes \Dh)$. More precisely, we show that all the elements of the Kasparov group $KK(\Dh, A\ot \Dh)$ are of the form
$[\varphi]-n[\iota]$ where $\varphi:\Dh \to \mathcal{K}\ot A\otimes \Dh$ is a
 $^*$-homomorphism and $\iota:\Dh \to  A\otimes \Dh$ is the
 inclusion $\iota(d)=\be_A\ot d$ and $n\in \mathbb{N}$.
 Moreover,  two non-zero $^*$-homomorphisms $\varphi,\psi: \Dh \to \mathcal{K}\ot A\ot \Dh$ with $\varphi(\be_\Dh)=\psi(\be_\Dh)=e$ have the same KK-theory class if and only
 if there is a unitary-valued continuous map $u:[0,1)\to e(\mathcal{K}\ot A\otimes \Dh)e$, $t\mapsto u_t$ such that $u_0=e$ and
 $\lim_{t\to 1}\|u_t\,\varphi(d)\,u_t^*-\psi(d)\|=0$ for all $d\in \Dh$.
In addition, we show that $KK_i(\Dh,\Dh\ot A)\cong K_i(\Dh\ot A)$,
$i=0,1$.

One may note the similarity to
the descriptions of $KK(\mathcal{O}_\infty,\mathcal{O}_\infty\ot A)$
(\cite{Kir:class},\cite{Phi:class}) and
 $KK(\mathbb{C},\mathbb{C}\ot A)$.
However, we do not require that $\Dh$ satisfies
the universal coefficient theorem (UCT) in KK-theory. 
In the same spirit, we characterize $\Oh_{2}$ and the universal UHF algebra $\Qh$ using $K$-theoretic conditions, but without involving the UCT. 
 
As another application of Theorem~\ref{D-asu} (and the results of \cite{HirshbergRordamWinter:absorb-ssa}), we prove in  \cite{Dadarlat-Winter:trivial-D-fields}  an automatic trivialization result for continuous fields
with strongly self-absorbing fibres over finite dimensional spaces.

The second named author would like to thank Eberhard Kirchberg for an inspiring conversation on the problem of proving Theorem~\ref{D-asu}.

\section{Strongly self-absorbing $C^{*}$-algebras}

In this section we recall the notion of strongly self-absorbing $C^{*}$-algebras and some facts from \cite{TomsWinter:ssa}.

\bn
\begin{dfs}
Let $A$, $B$ be $C^{*}$-algebras and $\sigma,\gamma: A \to B$ be $^{*}$-homomorphisms. Suppose that $B$ is unital.
\begin{enumerate}
\item We say that $\sigma$ and $\gamma$ are approximately unitarily equivalent, $\sigma \au \gamma$, if there is a sequence $(u_{n})_{n \in \N}$ of unitaries in  $B$ such that
\[
\|u_{n} \sigma(a) u_{n}^{*} - \gamma(a)\| \stackrel{n \to \infty}{\longrightarrow} 0
\]
for every $a \in A$.
If all $u_n$ can be chosen  to be in $\Uh_{0}(B)$, the connected
component of $\be_B$ of the unitary group $\Uh(B)$, then we say that
$\sigma$ and $\gamma$ are strongly approximately unitarily equivalent,
written $\sigma \sau \gamma$.
\item We say that $\sigma$ and $\gamma$ are asymptotically unitarily equivalent, $\sigma \asu \gamma$, if there is a norm-continuous path $(u_{t})_{t \in [0,\infty)}$ of unitaries in  $B$ such that
\[
\|u_{t} \sigma(a) u_{t}^{*} - \gamma(a)\| \stackrel{t \to \infty}{\longrightarrow} 0
\]
for every $a \in A$. If one can arrange that $u_0=\be_B$ and hence ($u_t\in\Uh_{0}(B)$ for all $t$), then we say that
$\sigma$ and $\gamma$ are strongly asymptotically unitarily equivalent,
written $\sigma \sasu \gamma$.
\end{enumerate}
\end{dfs}
\en

\bn
The concept of strongly self-absorbing $C^{*}$-algebras was formally introduced in \cite[Definition~1.3]{TomsWinter:ssa}:

\begin{dfs}
\label{def:ssa}
A separable unital $C^{*}$-algebra $\Dh$ is strongly self-absorbing, if $\Dh \neq \C$ and there is an isomorphism $\varphi: \Dh \to \Dh \otimes \Dh$ such that $\varphi \au \id_{\Dh}\otimes \be_{\Dh}$.
\end{dfs}
\en

\bn
Recall \cite[Corollary 1.12]{TomsWinter:ssa}:

\label{D-au}
\begin{props}
Let $A$ and $\Dh$ be  unital $C^{*}$-algebras, with $\Dh$  strongly self-absorbing.  Then, any two unital $^{*}$-homomorphisms $\sigma, \gamma:\Dh \to A \otimes \Dh$ are approximately unitarily equivalent. In particular, any two unital endomorphisms of $\Dh$ are approximately unitarily equivalent.
\end{props}

We note that the assumption that $A$ is separable which appears in the original statement of  \cite[Corollary 1.12]{TomsWinter:ssa} is  not necessary and  was not used in the proof.
\en

\bn \label{lemma:commutator}
\begin{lms} Let $\Dh$ be a  strongly
self-absorbing $C^{*}$-algebra. Then there is a sequence of unitaries
$(w_n)_{n \in \N}$ in the commutator subgroup of $\Uh(\Dh\ot\Dh)$ such that for all
$d\in \Dh$
$\|  w_n(d \ot \be_\Dh)w_n^*-\be_\Dh \otimes d\|\to 0$ as $n\to \infty$.
\end{lms}

\begin{nproof} Let $\Fh \subset \Dh$ be a finite normalized set and
let $\varepsilon>0$. By \cite[Prop.~1.5]{TomsWinter:ssa}
there is a unitary $u\in \Uh(\Dh \ot \Dh)$ such that
$\|u (d \otimes \be_{\Dh}) u^{*} - \be_{\Dh} \otimes d \| <\varepsilon$ for all $d\in\Fh$.
Let $\theta:\Dh\ot \Dh\to \Dh$ be a $^*$-isomorphism. Then
 $\|(\theta(u^*)\ot \be_{\Dh}) u (d \otimes \be_{\Dh}) u^{*}(\theta(u)\ot \be_{\Dh}) - \be_{\Dh} \otimes d \| <\varepsilon$ for all $d\in\Fh$.
 By Proposition~\ref{D-au} $\theta \ot \be_\Dh\au \id_{\Dh \ot \Dh}$
 and so there is a unitary $v\in \Uh(\Dh \ot \Dh)$ such that
 $\|\theta(u^*) \ot \be_\Dh-v u^* v^*\|<\varepsilon$ and hence
 $\|(\theta(u^*) \ot \be_\Dh)u-v u^* v^*u\|<\varepsilon$.
 Setting $w=v u^*v^*u$ we deduce that
 $\|w (d \otimes \be_{\Dh}) w^{*} - \be_{\Dh} \otimes d \| <3\varepsilon$ for all $d\in\Fh$.
\end{nproof}
\en

\bn
\label{connected-component-intertwining}
\begin{rems}
 In the situation of Proposition~\ref{D-au},
suppose that the commutator subgroup of $\Uh(\Dh)$ is contained in $\Uh_{0}(\Dh)$.
This will happen for instance if
$\Dh$ is  assumed to be $K_{1}$-injective. Then one may choose the unitaries  $(u_{n})_{n \in \N}$ which implement the approximate unitary equivalence between $\sigma$ and $\gamma$ to lie in $\Uh_{0}(A\ot \Dh)$. This follows from   \cite[(the proof of) Corollary 1.12]{TomsWinter:ssa}, since the unitaries  $(u_{n})_{n \in \N}$ are essentially images of the unitaries $(w_n)_{n \in \N}$ of Lemma~\ref{lemma:commutator} under suitable unital $^{*}$-homomorphisms.
\end{rems}
\en

\section{Asymptotic vs.\ approximate unitary equivalence}

It is the aim of this section to establish  a continuous version of Proposition \ref{D-au}.

\bn
\label{almost-central-homotopy}
\begin{lms}
Let $\Dh$ be separable unital strongly
self-absorbing $C^{*}$-algebra.  For any finite subset $\Fh \subset \Dh$ and
$\varepsilon>0$, there are a finite subset $\Gh \subset \Dh$ and $\delta>0$
such that the following holds:

If $A$ is another unital $C^{*}$-algebra and $\sigma:\Dh \to A \otimes \Dh$
is a unital $^{*}$-homomorphism, 
and if
 $w \in \Uh_0(A \otimes \Dh)$ is a unitary satisfying
\[
\|[w,\sigma(d)]\| < \delta
\]
for all $d \in \Gh$, then there is a continuous path $(w_{t})_{t \in [0,1]}$ of unitaries in $\Uh_0(A \otimes \Dh)$ such that $w_{0}=w$, $w_{1}=\be_{A \otimes \Dh}$ and
\[
\|[w_{t},\sigma(d)]\| < \varepsilon
\]
for all $d \in \Fh$, $t \in  [0,1]$.
\end{lms}

\begin{nproof}
We may clearly assume that the elements of $\Fh$ are  normalized
and that $\varepsilon <1$. Let $u \in \Dh \otimes \Dh$ be a unitary satisfying
\begin{equation}
\label{3}
\|u (d \otimes \be_{\Dh}) u^{*} - \be_{\Dh} \otimes d \| < \frac{\varepsilon}{20}
\end{equation}
for all $d \in \Fh$. There exist $k \in \N$ and elements $s_{1}, \ldots, s_{k},t_{1}, \ldots, t_{k} \in \Dh$ of norm at most one such that
\begin{equation}
\label{55}
\|u - \sum_{i=1}^{k} s_{i} \otimes t_{i} \| < \frac{\varepsilon}{20}.
\end{equation}
Set
\begin{equation}
\label{8}
\delta:= \frac{\varepsilon}{k \cdot 10}
\end{equation}
and
\begin{equation}
\label{7}
\Gh:= \{s_{1}, \ldots,s_{k}\} \subset \Dh.
\end{equation}
Now let $w \in \Uh_{0}(A \otimes \Dh)$ be a unitary as in the assertion of the lemma, i.e., $w$ satisfies
\begin{equation}
\label{6}
\|[w,\sigma(s_i)]\| < \delta
\end{equation}
for all $i=1,\dots,k$. We proceed to construct the path $(w_{t})_{t \in [0,1]}$.

By \cite[Remark 2.7]{TomsWinter:ssa} there is a unital $^{*}$-homomorphism
\[
\varphi: A \otimes \Dh \otimes \Dh \to A \otimes \Dh
\]
such that
\begin{equation}
\label{1}
\|\varphi(a \otimes \be_{\Dh}) - a\| < \frac{\varepsilon}{20}
\end{equation}
for all $a \in \sigma(\Fh) \cup \{w\}$.

Since $w \in \Uh_{0}(A \otimes \Dh)$, there is a path $(\bar{w}_{t})_{t \in [\halb,1]}$ of unitaries in $A \otimes \Dh$ such that
\begin{equation}
\label{4}
\bar{w}_{\halb}=w \mbox{ and } \bar{w}_{1} = \be_{A \otimes \Dh}.
\end{equation}
 For $t \in [\halb,1]$ define
\begin{equation}
\label{2}
w_{t}:= \varphi((\sigma \otimes \id_{\Dh})(u)^{*}(\bar{w}_{t} \otimes \be_{\Dh}) (\sigma \otimes \id_{\Dh})(u)) \in \Uh(A \otimes \Dh);
\end{equation}
then $(w_{t})_{t \in [\halb,1]}$ is a continuous path of unitaries in $A \otimes \Dh$. For $t \in [\halb,1]$ and $d \in \Fh$ we have
\begin{eqnarray}
\label{9}
\lefteqn{\|[w_{t},\sigma(d)]\|} \nonumber \\
& = & \|w_{t} \sigma(d) w_{t}^{*} - \sigma(d)\|   \nonumber \\
& \stackrel{\eqref{1}}{<} & \| w_{t} \varphi(\sigma(d) \otimes \be_{\Dh}) w_{t}^{*} - \varphi(\sigma(d) \otimes \be_{\Dh})\| + 2 \cdot \frac{\varepsilon}{20}  \nonumber \\
& \stackrel{\eqref{2}}{\le} & \| ((\sigma \otimes \id_{\Dh})(u))^{*}(\bar{w}_{t} \otimes \be_{\Dh}) ((\sigma \otimes \id_{\Dh})(u (d \otimes \be_{\Dh})u^{*})) (\bar{w}_{t}^{*} \otimes \be_{\Dh}) \nonumber \\
&& \cdot  ((\sigma \otimes \id_{\Dh})(u))  - ((\sigma \otimes \id_{\Dh})(d \otimes \be_{\Dh}))\| + \frac{\varepsilon}{10}  \nonumber \\
& \stackrel{\eqref{3}}{<} & \| ((\sigma \otimes \id_{\Dh})(u))^{*}(\bar{w}_{t} \otimes \be_{\Dh}) ((\sigma \otimes \id_{\Dh})(\be_{\Dh}\otimes d)) (\bar{w}_{t}^{*} \otimes \be_{\Dh}) \nonumber \\
&& \cdot ((\sigma \otimes \id_{\Dh})(u))  - ((\sigma\otimes \id_{\Dh})(d \otimes \be_{\Dh}))\| + \frac{\varepsilon}{10}+ \frac{\varepsilon}{20} \nonumber \\
& = & \| (\sigma \otimes \id_{\Dh})(u^{*} (\be_{\Dh}\otimes d)u - d \otimes \be_{\Dh})\| + \frac{\varepsilon}{10}+ \frac{\varepsilon}{20} \nonumber \\
& < & \frac{\varepsilon}{20} + \frac{\varepsilon}{10} + \frac{\varepsilon}{20}  \nonumber \\
& < & \frac{\varepsilon}{3},
\end{eqnarray}
where for the last equality we have used that the $\bar{w}_{t}$ are unitaries and that $\sigma$ is a unital $^{*}$-homomorphism. Furthermore, we have
\begin{eqnarray*}
\lefteqn{\|w_{\halb} - w\|} \\
& \stackrel{\eqref{4},\eqref{2}}{=} & \|\varphi(((\sigma \otimes \id_{\Dh})(u))^{*} (w \otimes \be_{\Dh}) ((\sigma \otimes \id_{\Dh})(u))) - w\| \\
& \stackrel{\eqref{55}}{<} & \| \varphi(((\sigma \otimes \id_{\Dh})(u))^{*} (w \otimes \be_{\Dh}) (\sum_{i=1}^{k} \sigma(s_{i}) \otimes t_{i})) - w\| + \frac{\varepsilon}{20} \\
& \le & \| \varphi(((\sigma \otimes \id_{\Dh})(u))^{*} (\sum_{i=1}^{k} \sigma(s_{i}) \otimes t_{i}) (w\otimes \be_{\Dh})) - w\| \\
&& + \sum_{i=1}^{k} \|[w, \sigma(s_{i})]\| \cdot \|t_{i}\| + \frac{\varepsilon}{20} \\
& \stackrel{\eqref{6},\eqref{7},\eqref{55}}{<} & \|\varphi (w \otimes \be_{\Dh}) - w\| + k \cdot \delta + 2 \cdot \frac{\varepsilon}{20} \\
& \stackrel{\eqref{1},\eqref{8}}{<} & \frac{\varepsilon}{20} + \frac{\varepsilon}{10} + 2 \cdot \frac{\varepsilon}{20} \\
& < & \frac{\varepsilon}{3} .
\end{eqnarray*}
The above estimate allows us to extend the path $(w_{t})_{t \in [\halb,1]}$ to the whole interval $[0,1]$ in the desired way: We have $\|w_{\halb} w^{*} - \be_{\Dh}\| < \frac{\varepsilon}{3} <2$, whence
$-1$ is not in the spectrum of $w_{\halb} w^{*}$. By functional calculus, there is $a=a^*\in A \otimes \Dh$ with $\|a\|<1$ such that
$w_{\halb} w^{*}=\exp(\pi ia)$.
 For $t \in [0,\halb)$ we may therefore define a continuous path of unitaries
\[
w_{t}:= (\exp(2\pi i t a)) w \in \Uh(A \otimes \Dh).
\]
It is clear that $w_{0}=w$ and $w_{t} \to w_{\halb}$ as $t \to (\halb)_{-}$, whence $(w_{t})_{t \in [0,1]}$ is a continuous path of unitaries in $A$ satisfying $w_{0} = w$ and $w_1=\be_A \otimes \Dh$. Moreover, it is easy to see that
\[
\|w_{t} - w\| \le \|w_{\halb} - w\| < \frac{\varepsilon}{3}
\]
for all $t \in [0,\halb)$, whence
\[
\|[w_{t} ,\sigma(d)]\| < \|[w_{\halb},\sigma(d)]\| + \frac{2}{3} \,\varepsilon \stackrel{\eqref{9}}{<} \varepsilon
\]
for $t \in [0, \halb)$, $d \in \Fh$.

We have now constructed a path $(w_{t})_{t \in [0,1]} \subset \Uh(A)$  with the desired properties.
\end{nproof}
\en

\bn
\label{D-asu}
\begin{thms}
Let $A$ and $\Dh$ be  unital $C^{*}$-algebras,
with $\Dh$ separable, strongly self-absorbing and $K_{1}$-injective.
 Then, any two unital $^{*}$-homomorphisms $\sigma, \gamma:\Dh \to A \otimes \Dh$ are strongly asymptotically unitarily equivalent.  In particular, any two unital endomorphisms of $\Dh$ are strongly asymptotically unitarily equivalent.
\end{thms}

\begin{nproof}
Note that the second statement follows from the first one with $A =\Dh$, since $\Dh \cong \Dh \otimes \Dh$ by assumption.

Let $A$ be a unital $C^{*}$-algebra such that $A\cong A\ot \Dh$ and
 let $\sigma,\gamma:\Dh \to A$ be unital $^*$-homomorphisms.
 We shall prove that $\sigma$ and $\gamma$ are strongly asymptotically unitarily equivalent.
 Choose an increasing sequence
\[
\Fh_{0} \subset \Fh_{1} \subset \ldots
\]
of finite subsets of $\Dh$ such that $\bigcup \Fh_{n}$ is a dense subset of  $\Dh$. Let $1> \varepsilon_{0} > \varepsilon_{1} > \ldots $ be a decreasing sequence of strictly positive numbers converging to $0$.

For each $n \in \N$, employ Lemma \ref{almost-central-homotopy} (with $\Fh_{n}$ and $\varepsilon_{n}$ in place of $\Fh$ and $\varepsilon$) to obtain a finite subset $\Gh_{n} \subset \Dh$ and $\delta_{n}>0$. We may clearly assume that
\begin{equation}
\label{G_n-delta_n}
\Fh_{n} \subset \Gh_{n} \subset \Gh_{n+1} \mbox{ and that } \delta_{n+1} < \delta_{n} < \varepsilon_{n}
\end{equation}
for all $n \in \N$.

Since $\sigma$ and $\gamma$ are strongly approximately unitarily equivalent by Proposition \ref{D-au} and Remark \ref{connected-component-intertwining}, there is a sequence of unitaries $(u_{n})_{n \in \N} \subset \Uh_0(A)$ such that
\begin{equation}
\label{gamma-sigma-au}
\|u_{n} \sigma(d) u_{n}^{*} - \gamma(d)\| < \frac{\delta_{n}}{2}
\end{equation}
for all $d \in \Gh_{n}$ and $n \in \N$.
Let us set
\[
w_{n}:= u_{n+1}^{*}u_{n}, \, n \in \N.
\]
Then $w_{n} \in \Uh_{0}(A)$ and
\begin{eqnarray*}
\lefteqn{\|[w_{n},\sigma(d)]\|}\\
& = & \|w_{n} \sigma(d) w_{n}^{*} - \sigma(d)\| \\
& \le & \|u_{n+1}^{*} u_{n} \sigma(d) u_{n}^{*}u_{n+1} - u_{n+1}^{*} \gamma(d) u_{n+1}\| \\
& & + \|u_{n+1}^{*} \gamma(d) u_{n+1} - \sigma(d) \| \\
& < & \frac{\delta_{n}}{2} + \frac{\delta_{n+1}}{2} \\
& < & \delta_{n}
\end{eqnarray*}
for $d \in \Gh_{n}$, $n \in \N$. Now by Lemma \ref{almost-central-homotopy} (and the choice of the $\Gh_{n}$ and $\delta_{n}$), for each $n$ there is a continuous path $(w_{n,t})_{t \in [0,1]}$ of unitaries in $\Uh_0(A)$ such that $w_{n,0}=w_{n}$, $w_{n,1}=\be_{A}$ and
\begin{equation}
\label{almost-central-unitaries}
\|[w_{n,t}, \sigma(d)]\| < \varepsilon_{n}
\end{equation}
for all $d \in \Fh_{n}$, $t \in [0,1]$.

Next, define a path $(\bar{u}_{t})_{t \in [0,\infty)}$ of unitaries in $\Uh_0(A)$ by
\[
\bar{u}_{t}:= u_{n+1}w_{n,t-n} \mbox{ if } t \in [n,n+1).
\]
We have that
\begin{equation}
\label{10}
\bar{u}_{n}=u_{n+1} w_{n} = u_{n}
\end{equation}
and that
\[
\bar{u}_{t} \to u_{n+1}
\]
as $t \to n+1$ from below, which implies that the path $(\bar{u}_{t})_{t \in [0,\infty)}$ is continuous in $\Uh_0(A)$. Furthermore, for $t \in [n,n+1)$ and $d \in \Fh_{n}$ we obtain
\begin{eqnarray*}
\lefteqn{\|\bar{u}_{t} \sigma(d) \bar{u}_{t}^{*} - \gamma(d)\|}\\
& = & \|u_{n+1} w_{n,t-n} \sigma(d) w_{n,t-n}^{*} u_{n+1}^{*} - \gamma(d)\| \\
& \stackrel{\eqref{almost-central-unitaries}}{<} & \|u_{n+1} \sigma(d) u_{n+1}^{*} - \gamma(d)\| + \varepsilon_{n} \\
& \stackrel{\eqref{gamma-sigma-au},\eqref{G_n-delta_n}}{<} & \frac{\delta_{n+1}}{2} + \varepsilon_{n}\\
& \stackrel{\eqref{G_n-delta_n}}{<} & 2 \varepsilon_{n}.
\end{eqnarray*}
Since the $\Fh_{n}$ are nested and the $\varepsilon_{n}$ converge to $0$, we have
\begin{equation}
\label{gamma-sigma-asu}
\|\bar{u}_{t} \sigma(d) \bar{u}_{t}^{*} - \gamma(d)\| \stackrel{t \to \infty}{\longrightarrow} 0
\end{equation}
for all $d \in \bigcup_{n=0}^{\infty} \Fh_{n}$; by continuity and since $\bigcup_{n=0}^{\infty} \Fh_{n}$ is dense in  $\Dh$, we have \eqref{gamma-sigma-asu} for all $d \in \Dh$.
Since $\bar{u}_{0}\in \Uh_0(A)$ we may arrange that $\bar{u}_{0}=\be_A$.
\end{nproof}
\en

\section{The group $KK(\Dh,A\ot \Dh)$ and  some applications}

\bn
For a separable $C^{*}$-algebra $\Dh$ we endow the group of automorphisms $\Aut(\Dh)$ with
the point-norm topology.

\begin{cors}\label{Cor:aut=contractible}
    Let $\Dh$ be a separable, unital, strongly
self-absorbing and $K_{1}$-injective $C^{*}$-algebra.
Then
$[X,\mathrm{Aut}(\Dh)]$ reduces to a point for any compact Hausdorff space
$X$.
\end{cors}

\begin{nproof} Let $\varphi,\psi:X \to \Aut(\Dh)$ be continuous
maps. We identify $\varphi$ and $\psi$ with unital $^*$-homomorphisms
$\varphi,\psi:\Dh\to \Ch(X)\ot \Dh$. By Theorem~\ref{D-asu}, $\varphi$ is
strongly asymptotically unitarily equivalent to $\psi$. This gives a homotopy between the
two maps $\varphi,\psi:X \to \Aut(\Dh)$.
\end{nproof}
\en

\bn
\begin{rems}
The conclusion of Corollary~\ref{Cor:aut=contractible}
was known before for $\Dh$ a UHF algebra of infinite type and $X$ a CW complex by \cite{Thomsen:homotopy-aut-UHF}, for
$\Dh=\mathcal{O}_2$ by \cite{Kir:class} and \cite{Phi:class}, and
for $\Dh=\mathcal{O}_\infty$ by \cite{Dad:bundles-fdspaces}. It is new for the Jiang--Su algebra. 
\end{rems}
\en

\bn
For unital $C^{*}$-algebras $\Dh$ and $B$ we denote by $[\Dh,B]$ the set of  homotopy classes of
 unital $^*$-homomorphisms from $\Dh$ to $B$. By a similar argument as above we also have the following corollary.

\begin{cors}\label{Cor:homotopy=singleton}
    Let $\Dh$ and $A$ be  unital $C^{*}$-algebras.
If $\Dh$ is separable, strongly self-absorbing and $K_{1}$-injective, then
$[\Dh,A\otimes \Dh]$ reduces to a singleton.
\end{cors}
\en

\bn
For separable unital $C^{*}$-algebras $\Dh$ and $B$, let $\chi_i:KK_i(\Dh,B)\to
KK_i(\mathbb{C},B)\cong K_i(B)$, $i=0,1$  be the morphism of groups
induced by the unital inclusion $\nu:\mathbb{C}\to \Dh$.

\begin{thms}\label{Thm:KK-topology}
    Let $\Dh$  be a unital, separable
and strongly self-absorbing $C^{*}$-algebra. Then for any separable $C^{*}$-algebra $A$, the map
 $\chi_i:KK_i(\Dh,A\otimes \Dh)\to K_i(A\otimes \Dh)$
is bijective, for $i=0,1$. In particular both groups  $KK_i(\Dh,A\otimes D) $ are
countable and discrete with respect to their  natural topology.
\end{thms}

\begin{nproof} Since $\Dh$ is KK-equivalent to $\Dh\otimes \mathcal{O}_\infty$, we may assume
 that $\Dh$ is  purely infinite and in particular $K_1$-injective by
\cite[Prop.~4.1.4]{Ror:encyclopedia}.
 Let $C_\nu \Dh$ denote the mapping cone $C^*$-algebra of $\nu$.
 By \cite[Cor.~3.10]{Dad:homotopy-aut},  there is a
bijection $[\Dh,A\otimes \Dh]\to KK(C_\nu \Dh, SA\otimes \Dh)$ and hence $KK(C_\nu \Dh,
SA\otimes \Dh)=0$ for all separable and unital $C^*$-algebras $A$ as a consequence of
Corollary~\ref{Cor:homotopy=singleton}. Since  $KK(C_\nu \Dh, A\otimes \Dh)$ is
isomorphic to $KK(C_\nu \Dh, S^2A\otimes \Dh)$ by Bott periodicity and the latter
group injects in $KK(C_\nu \Dh, SC(\mathbb{T})\otimes A\otimes \Dh)=0$, we have that
$KK_i(C_\nu \Dh,\Dh\otimes A)=0$ for all unital and separable $C^*$-algebras $A$ and
$i=0,1$. Since $KK_i(C_\nu \Dh,\Dh\otimes A)$ is a subgroup
of $KK_i(C_\nu \Dh,\Dh\otimes \widetilde{A})=0$ (where $\widetilde{A}$ is the unitization of $A$) we see that
$KK_i(C_\nu \Dh,\Dh\otimes A)=0$ for all separable $C^*$-algebras $A$.
 Using the Puppe exact sequence, where $\chi_i=\nu^*$,
\[ \xymatrix@C=1.2em{
 KK_{i+1}(C_\nu \Dh,A\otimes \Dh)\ar[r]
                & KK_i(\Dh,A\otimes \Dh)
                \ar[r]^{\chi_i}&KK_i(\mathbb{C},A\otimes \Dh)\ar[r]
                &{KK_{i}(C_\nu \Dh,A\otimes \Dh)}
 }\]
we conclude that  $\chi_i$ is an isomorphism, $i=0,1$.
 The map $\chi_i=\nu^*$ is continuous
 since it is given by the Kasparov product with a fixed element
(we refer the reader   to
 \cite{Sc:fine1}, \cite{Pim:KK-top} or
 \cite{Dad:kk-top} for a background on the topology of the Kasparov groups).
 Since the topology of $K_i$ is discrete and $\chi_i$ is injective, it
 follows that the topology of $KK_i(\Dh,A\ot D)$ is also discrete.
 The countability of $KK_i(\Dh,A\ot D)$ follows from that of $K_i(A\ot D)$,
 as $A\otimes \Dh$ is separable.
\end{nproof}
\en

\bn
\begin{rems}
In contrast to Theorem~\ref{Thm:KK-topology}, if $\Dh$ is the universal UHF algebra, then
$KK(\Dh,\mathbb{C})\cong \mathrm{Ext}(\mathbb{Q},\mathbb{Z})\cong \mathbb{Q}^\mathbb{N}$
has the power of the continuum \cite[p.~221]{Fuc:inf}.
\end{rems}
\en

\bn
Let $\Dh$ and $A$ be as in Theorem~\ref{Thm:KK-topology} and assume in addition that $\Dh$ is $K_1$-injective and  $A$ is unital.
Let $\iota:\Dh \to A\ot \Dh$ be defined by
$\iota(d)=\be_A \ot d$.

\begin{cors} If $e\in\mathcal{K}\ot A\ot \Dh$ is a projection,
and $\varphi,\psi:\Dh \to e(\mathcal{K}\ot A\ot \Dh)e$ are two unital $^*$-homomorphisms,
then $\varphi \sasu \psi$ and hence $[\varphi]=[\psi]\in KK(\Dh,A\otimes \Dh)$. Moreover:
 \[KK(\Dh,A\otimes \Dh)=\{[\varphi]-n[\iota]\,|\, \varphi:\Dh \to \mathcal{K}\ot A\ot \Dh \,\mbox{is a $^*$-homomorphism}, \,n \in \mathbb{N}\}.\]
\end{cors}

\begin{nproof} Let $\varphi$, $\psi$  and $e$ be as in the first part of the statement. By \cite[Cor.~3.1]{TomsWinter:ssa}, the unital $C^{*}$-algebra $e(\mathcal{K}\ot A\ot \Dh)e$ is $\Dh$-stable, being a hereditary subalgebra of a $\Dh$-stable $C^*$-algebra.
Therefore $\varphi \sasu \psi$ by Theorem~\ref{D-asu}.

Now for the second part of the statement, let $x\in KK(\Dh,A\otimes \Dh)$ be an arbitrary element. Then $\chi_0(x)=[e]-n[\be_{A\otimes \Dh}]$
for some projection $e\in\mathcal{K}\ot A\ot \Dh$ and $n \in \mathbb{N}$. Since $e(\mathcal{K}\ot A\ot \Dh)e$ is $\Dh$-stable,
there is a unital $^*$-homomorphism $\varphi:\Dh \to e(\mathcal{K}\ot A\ot \Dh)e$. Then
\[\chi_0([\varphi]-n[\iota])=[\varphi(\be_\Dh)]-n[\iota(\be_\Dh)]=
[e]-n[\be_{A\otimes \Dh}]=\chi_0(x),\]
and hence $[\varphi]-n[\iota]=x$ since $\chi_0$ is injective
by Theorem~\ref{Thm:KK-topology}.
\end{nproof}
\en

In the remainder of the paper we give characterizations
for the Cuntz algebra $\mathcal{O}_2$ and for the universal UHF-algebra which do not require the UCT. The latter result is a variation
of a theorem of Effros and Rosenberg \cite{Effros-Rosen:flip}.

\bn
\begin{props}\label{Prop:O(2)}
    Let $\Dh$ be a separable unital strongly self-absorbing $C^{*}$-algebra.
If $[\be_{\Dh}]=0$ in $K_0(\Dh)$, then $\Dh\cong \mathcal{O}_2$.
\end{props}

\begin{nproof} Since $\Dh$ must be nuclear (see \cite{TomsWinter:ssa}),
$\Dh$ embeds unitally in $\mathcal{O}_2$ by Kirchberg's theorem.
 $\Dh$ is not stably finite since $[\be_{\Dh}]=0$. By
the dichotomy of \cite[Thm.~1.7]{TomsWinter:ssa} $\Dh$ must be purely infinite.
Since $[\be_{\Dh}]=0$ in $K_0(\Dh)$, there is a unital embedding $\mathcal{O}_2\to \Dh$,
see \cite[Prop.~4.2.3]{Ror:encyclopedia}.
 We
conclude that $\Dh$ is isomorphic to $\mathcal{O}_2$ by
\cite[Prop.~5.12]{TomsWinter:ssa}.
\end{nproof}
\en

\bn
\begin{props} \label{approx-absorb}
Let $\Dh$, $A$ be separable, unital, strongly self-absorbing
$C^*$-algebras. Suppose that for any finite subset $\mathcal{F}$ of $\Dh$ and any
$\varepsilon>0$ there is a u.c.p. map $\varphi:\Dh\to A$ such that $\|\varphi(c
d)-\varphi(c)\varphi(d)\|<\varepsilon$ for all $c,d\in \mathcal{F}$. Then $A\cong
A\otimes \Dh$.
\end{props}

\begin{nproof}
By \cite[Thm.~2.2]{TomsWinter:ssa} it suffices to show that for any given finite
subsets $\mathcal{F}$ of $\Dh$, $\mathcal{G}$ of $A$ and any $\varepsilon>0$ there
is u.c.p. map $\Phi:\Dh\to A$ such that (i) $\|\Phi(c d)-\Phi(c)\Phi(d)\|<\varepsilon$
for all $c,d\in \mathcal{F}$ and (ii) $\|[\Phi(d),a]\|<\varepsilon$ for all $d\in
\mathcal{F}$ and $a \in \mathcal{G}$. We may  assume that $\|d\|\leq 1$ for
all $d\in \mathcal{F}$.
Since $A$ is strongly self-absorbing, by \cite[Prop.~1.10]{TomsWinter:ssa} there
is a unital $^*$-homomorphism $\gamma:A\otimes A \to A$ such that $\|\gamma(a\ot
\be_A)-a\|<\varepsilon/2$ for all $a \in \mathcal{G}$. On the other hand, by
assumption there is a u.c.p. map $\varphi:\Dh\to A$ such that $\|\varphi(c
d)-\varphi(c)\varphi(d)\|<\varepsilon$ for all $c,d\in \mathcal{F}$. Let us
define a u.c.p. map $\Phi:\Dh\to A$ by $\Phi(d)=\gamma(\be_A\otimes \varphi(d))$. It is
clear that $\Phi$ satisfies (i) since $\gamma$ is a $^*$-homomorphism. To conclude
the proof we check now that $\Phi$ also satisfies (ii). Let $d\in \mathcal{F}$
and $a \in \mathcal{G}$. Then
\begin{eqnarray*}
\lefteqn{\|[\Phi(d),a]\|}\\
& \leq  & \|[\Phi(d),a-\gamma(a\otimes \be_A)]\|+\|[\Phi(d),\gamma(a\otimes \be_A)]\|\\
& \leq & 2\|\Phi(d)\|\|a-\gamma(a\otimes \be_A)\| + \|[\gamma(\be_A\otimes \varphi(d)),\gamma(a\otimes \be_A)]\| \\
& < & 2 \varepsilon/2+0= \varepsilon.
\end{eqnarray*}
\end{nproof}
\en

\bn
\begin{props}\label{Prop:qd-ssa} Let $\Dh$ be a separable, unital, strongly self-absorbing
$C^*$-algebra. Suppose that $\Dh$ is quasidiagonal, it has cancellation of projections
and that $[\be_{\Dh}]\in n K_0(\Dh)^+$ for all $n \geq 1$. Then $\Dh$ is isomorphic to the
universal UHF algebra $\Qh$ with $K_0(\Qh)\cong \mathbb{Q}$.
\end{props}

\begin{nproof} Since $\Dh$ is separable unital and quasidiagonal, there is a unital $^*$-representation
 $\pi:\Dh \to B(H)$ on a separable Hilbert space $H$ and a sequence of nonzero
  projections $p_n\in B(H)$  of finite rank  $k(n)$ such that
  $\lim_{n\to \infty}\|[p_n,\pi(d)]\|=0$
 for all $d\in \Dh$. Then the sequence of u.c.p. maps
  $\varphi_n:\Dh \to p_n B(H)p_n\cong M_{k(n)}(\mathbb{C})\subset \Qh$
  is asymptotically multiplicative, i.e $\lim_{n\to \infty}\|
  \varphi_n(cd)-\varphi_n(c)\varphi_n(d))\|=0$
 for all $c,d\in \Dh$. Therefore $\Qh\cong \Qh\otimes \Dh$ by
 Proposition~\ref{approx-absorb}.

In the second part of the proof we show that $\Dh\cong \Dh\otimes \Qh$. Let $E_n:\Qh\to
M_{n!}(\mathbb{C})\subset \Qh$ be a conditional expectation onto
$M_{n!}(\mathbb{C})$. Then $\lim_{n\to \infty}\|E_n(a)-a\|=0$ for all $a\in \Qh$.

 By assumption, for each $n$ there is a projection $e$ in  $\Dh\otimes M_m(\mathbb{C})$ (for some $m$)
 such that $n![e]=[\be_{\Dh}]$ in $K_0(\Dh)$. Let $\varphi:M_{n!}(\mathbb{C})\to
 M_{n!}(\mathbb{C})\otimes e(\Dh\otimes M_m(\mathbb{C} ))e$ be defined by
 $\varphi(b)=b\otimes e$. Since $\Dh$ has cancellation of projections
 and since $n! [e]=[\be_{\Dh}]$, there is a partial isometry $v \in M_{n!}(\mathbb{C}) \otimes D\otimes
 M_m(\mathbb{C}) $ such that $v^*v=\be_{M_{n!}(\mathbb{C})}\otimes e$  and $vv^*=e_{11}\otimes \be_{\Dh}\otimes e_{11}$.
 Therefore $b \mapsto v\,\varphi(b)\,v^*$ gives a unital embedding of
 $M_{n!}(\mathbb{C})$ into $\Dh$. Finally, $\psi_n(a)=v\,(\varphi \circ E_n(a))\,v^*$
 defines a sequence of asymptotically multiplicative u.c.p. maps $\Qh\to \Dh$.
Therefore $\Dh\cong \Dh\otimes \Qh$ by
 Proposition~\ref{approx-absorb}.
\end{nproof}
\en

\bn
\begin{rems} Let $\Dh$ be a separable, unital, strongly self-absorbing
 and quasidiagonal $C^*$-algebra. Then $\Dh\otimes \Qh \cong \Qh$ by the first part
 of the proof of Proposition \ref{Prop:qd-ssa}. In particular
 $K_1(\Dh)\otimes \mathbb{Q} =0$ and
 $K_0(\Dh)\otimes \mathbb{Q}\cong \mathbb{Q}$ by the K\"{u}nneth formula
 (or by writing $\Qh$ as an inductive limit of matrices).
\end{rems}
\en

 \providecommand{\bysame}{\leavevmode\hbox to3em{\hrulefill}\thinspace}
\providecommand{\MR}{\relax\ifhmode\unskip\space\fi MR }
\providecommand{\MRhref}[2]{%
  \href{http://www.ams.org/mathscinet-getitem?mr=#1}{#2}
}
\providecommand{\href}[2]{#2}

%\bibliographystyle{amsplain}
%\bibliography{operator}

\end{document}